%% file: CDC14.tex
\documentclass[letterpaper, 12 pt, conference]{ieeeconf}

\IEEEoverridecommandlockouts     
\usepackage{graphicx}
\usepackage{multirow}
\usepackage{epstopdf}

\usepackage[cmex10]{amsmath}
\usepackage{amsfonts}
\usepackage{amssymb}

\usepackage{array}
\usepackage{color}
\usepackage{float}

\usepackage{graphicx,multirow,enumerate,etoolbox}
\usepackage{amsmath,amsfonts,amssymb,array,subfigure,color, algorithm, algorithmic}

\AtBeginEnvironment{align}{\setcounter{subeqn}{0}}
\newcounter{subeqn} \renewcommand{\thesubeqn}{\theequation\alph{subeqn}}%
\newcommand{\subeqn}{%
  \refstepcounter{subeqn}
  \tag{\thesubeqn}
}

\input defs.tex

\newtheorem{Theo}{Theorem}
\newtheorem{Problem}{Problem}
\newtheorem{Exa}{\textbf{Example}}

\def\norm#1{\|#1\|}

\newcommand{\seq}{\reals^{\mathbb{N}}}

\title{
Reconstruction of Support of a Measure From Its Moments
}

\author{A. M. Jasour and C. Lagoa
\thanks{}
\thanks{A. M Jasour is with the Department of Electrical Engineering,
        The Pennsylvania State University, UP, PA, 16802 
        {\tt\small jasour@psu.edu}}%
\thanks{C. Lagoa is with Faculty of Electrical Engineering, 
        The Pennsylvania State University, UP, PA, 16802 
        {\tt\small lagoa@psu.edu}}%
\thanks{
 }%
}

\begin{document}
\onecolumn
\maketitle
\thispagestyle{empty}
\pagestyle{empty}


\begin{abstract}

In this paper, we address the problem of reconstruction of support of a measure from its moments. More precisely, given a finite subset of  the moments of a measure, we develop a semidefinite program for approximating the support of measure using level sets of polynomials. To solve this problem, a sequence of convex relaxations is provided, whose optimal solution is shown to converge to the support of measure of interest. 
{Moreover, the  provided approach is modified to improve the results for uniform measures.}
Numerical examples are presented to illustrate the performance of the proposed approach.

\end{abstract}


\section{Introduction} \label{Sec_Intro}
In this paper, we aim at solving the problem of reconstructing of support of a measure using only its moments. More precisely, we consider the following problem.\\

\begin{Problem} \label{Prob1}
Given the moment sequence of a measure $\mu$, {
find a polynomial $\mathcal{P}: \mathbb{R}^n \rightarrow \mathbb{R}$ such that the set
\[
\cK = \{ x\in \mathbb{R}^n: \mathcal{P}(x)\geq 1 \}
\]
coincides with the support set of the measure $\mu$. \\
}
\end{Problem}

This problem has many applications in many different areas. A few examples are, problem of shape reconstruction from indirect measurements (\cite{ref_Int_shape2, ref_Int_shape1}), signal reconstruction from sparse measurements (\cite{ref_Int_signal2, ref_Int_signal1}), and problems in statistics \cite{ref_Int_stat1}.
Moreover, this problem can be applied in area of optimization. For example, moment approaches to polynomial optimization over semialgebraic sets where one aims at solving  
 \begin{equation}\label{M1}
\mathbf{f^*} =\sup_{\rm x \in \mathbb{K}} f(x) 
\end{equation} 
by looking at the moments of the measures in the following problem
 \begin{align}
\mathbf{\rho ^*} =\sup_{\rm \mu_x \in \cM(\mathbb{K})} \int_{\mathbb{K}} f(x) d\mu_x \\
\hbox{s.t.}\quad & \int_{\mathbb{K}} d\mu_x =1 \subeqn
 \end{align}
requires one to extract an  optimal solution by finding an $x^*$ in the support set of the optimal solution $\mu_x^*$ of the problem above; see \cite{ref_Int_opt1, ref_Int_opt2}. The same problem appears in semi-algebraic chance optimization problems of the form
\begin{equation} \label{M2}
\mathbf{P^*} =\sup_{\rm x\in {\chi}} Prob_{\mu_q} \{ q\in\reals^m:\  f(x,q)\geq \gamma \big\}
\end{equation}
which can also be solved using a moment approach and also require finding a point in the support set of a measure of which one only knows a finite set of moments; see \cite{ref_Int_chance2, ref_Int_chance1, ref_Int_chance3}.

In this paper, to reconstruct the support of the measure of interest from its moments, we develop a sequence of semidefinite programming (SDP) problems whose solutions converge to the solution of Problem \ref{Prob1}.

Several approaches have been proposed to construct the support from the moments information. In \cite{ref_Pre_1} an approach to exact reconstruction of convex polytope supports is proposed, which is based on the collection of moment formulas combined with Vandermonde factorization of finite rank Hankel matrices. In \cite{ref_Pre_2}, a method to reconstruct planar semi-analytic domains from their moments is proposed based on the diagonal Pade approximation where it can approximate arbitrarily closely any bounded domain. \cite{ref_Pre_3} provides an method to obtain a polynomial that vanishes on the boundary of support.

In this paper, we take a different approach. The proposed method relies on results on Sum of Squares (SOS) polynomials and also, results on necessary and sufficient condition for moment sequence to have a representing measure. A hierarchy of semidefinite relaxations for approximation of the support set is proposed. 

The outline of the paper is as follows. In Section \ref{Sec_Notation}, the notation used in this paper as well as preliminary results on measures theory and SOS polynomials are presented. In Section \ref{Sec_Convex}, a convex formulation of support reconstruction problem as well as numerical examples is provided . In Section \ref{Sec_Uniform} a modified SDP is given to improve the results for uniform measures. Concluding remarks are provided in Section \ref{Sec_Con}.


\section{ Notation and Preliminary Results} \label{Sec_Notation}

\subsection{ Notations and Definitions}

\label{sec:definitions}
Let $\mathbb{R}[x]$ be the ring of real polynomials in the variables $x \in \mathbb{R}^n$. 
Given $\cP\in\mathbb{R}[x]$, we will represent $\cP$ as $\sum_{\alpha\in\mathbb{N}^n} p_\alpha x^\alpha$ using the standard basis $\{x^\alpha\}_{\alpha\in \mathbb{N}^n}$ of $\mathbb{R}[x]$, and $\mathbf{p}=\{p_\alpha\}_{\alpha\in\mathbb{N}^n}$ denotes the sequence of polynomial coefficients. Moreover, let $ \Sigma^2[x] \subset \mathbb R [x]$ be the set of sum of squares (SOS) polynomials. $\sigma:\reals^n\rightarrow\reals$ is a SOS polynomial if it can be written as a sum of \emph{finitely} many squared polynomials, i.e. $\sigma (x)= \sum_{j=1}^{\ell} h_j(x)^2$ for some $\ell<\infty$ and $h_j\in\reals[x]$ for $1\leq j\leq \ell$. Given $n$ and $r$ in $\mathbb{N}$, we define $S_{n,r} := \binom{r+n}{n}$ and $\mathbb{N} ^{\rm n}_r = \{\alpha \in \mathbb N^n : \norm{\alpha}_1 \leq r \}$.
Let $\mathbb R_{\rm d}[x] \subset \mathbb R [x]$ denote the set of polynomials of degree at most $d\in \mathbb{N}$, which is indeed a vector space of dimension $S_{n,d}$. 

Let $\seq$ denote the space of real sequences, and let $\cM(\cK)$ be the set of finite Borel measures $\mu$ such that $supp(\mu)\subset\cK$, where $supp(\mu)$ denotes the support of the measure $\mu$; i.e., the smallest set that contains all measurable sets with strictly positive $\mu$ measure. A sequence $\mathbf y = \{ y_ \alpha \}_{\alpha\in\mathbb{N}^n}\in\seq$
is said to have a \emph{representing measure}, if there exists a finite Borel measure $\mu$ on $\reals^n$ such that $y_{\alpha } = \int{x^{\alpha} d\mu }$ for every $\alpha \in \mathbb N ^n$ -see~\cite{ref_Int_opt1, ref_Int_opt2}. In this case, $ \mathbf y $ is called the moment sequence of the measure $\mu $. Given a square symmetric matrices $A$, the notation $A \succcurlyeq 0$ denotes positive semidefiniteness of $A$.

\textbf{Putinar’s property:} A closed semialgebraic set $\cK = \{ x\in \mathbb{R}^n: \mathcal{P}_j(x)\geq0,\ j=1,2,\dots ,\ell\ \}$ defined by polynomials $\mathcal{P}_j\in \mathbb R [x]$ satisfies \emph{Putinar's property} if there exists $\mathcal{U}\in \mathbb R [x]$ such that $\lbrace x:  \mathcal{U}(x) \geq 0 \rbrace $ is compact and $\mathcal{U} = \mathcal{\sigma}_0 + \sum_{j=1}^{m} \mathcal{\sigma}_j\mathcal{P}_j $ for some SOS polynomials $\sigma_j \subset \Sigma^2[x]$ -- see~\cite{ref_PreR_1, ref_Int_opt2}. Putinar's property holds if the level set $\lbrace x: \mathcal{P}_j(x) \geq 0 \rbrace$ is compact for some $j$, or if all $ \mathcal{P}_j $
are affine and $\cK $ is compact - see~\cite{ref_PreR_1}. Clearly these results imply that if there exits $M>0$ such that the polynomial $\cP_{\ell+1}(x):= M- \Vert x \Vert^2 \geq 0 $ for all $x\in\cK$, then $\cK\cap\{x: \cP_{\ell+1}\geq 0\}$ satisfies Putinar's property.

\textbf{Moment matrix:} Given $r\geq 1$ and the sequence $\{y_\alpha\}_{\alpha\in\mathbb{N}^n}$, the moment matrix $M_r({\mathbf y})\in\reals^{S_{n,r}\times S_{n,r}}$, containing all the moments up to order $2r$, is a symmetric matrix defined as follows~\cite{ref_Int_opt1, ref_Int_opt2}:
\begin{equation}\label{momnt matirx def}
M_r ( \mathbf y )(i,j)=y_{\alpha^{(i)}+\alpha^{(j)}},\ \ \ 1 \leq i,j \leq S_{n,r},
\end{equation}
where the elements of the moment sequence $ \mathbf y =\{ y_\alpha\}_{\alpha\in\mathbb{N}^n}$ are sorted according to a graded reverse lexicographic order of the corresponding monomials so that we have $\reals^n\ni\mathbf{0} = \alpha ^{(1)} < \ldots < \alpha ^{(S_{n,2r})}$ and $S_{n,2r}$ is the number of moments in $\mathbb{R}^n$ up to order $2r$.

For $r=2$ and $n=2$, the moment matrix containing moments up to order $2r$ is
\begin{equation} \label{moment matrix exa}
M_2\left({\mathbf y}\right)=\left[ \begin{array}{c}

\begin{array}{ccc} y_{00} \ | & y_{10} & y_{01}| \end{array}
\begin{array}{ccc} y_{20} & y_{11} & y_{02} \end{array}
 \\
 \begin{array}{ccc} - & - & - \end{array}
\ \ \ \  \begin{array}{ccc} - & - & - \end{array}
 \\

 \begin{array}{ccc} y_{10}\ | & y_{20} & y_{11}| \end{array}
\ \begin{array}{ccc} y_{30} & y_{21} & y_{12} \end{array}
 \\

 \begin{array}{ccc} y_{01}\ | & y_{11} & y_{02}| \end{array}
\ \begin{array}{ccc} y_{21} & y_{12} & y_{03} \end{array}
 \\

 \begin{array}{ccc} - & - & - \end{array}
\ \ \ \ \  \begin{array}{ccc} - & - & - \end{array}
 \\

 \begin{array}{ccc} y_{20}\ | & y_{30} & y_{21}| \end{array}
\ \begin{array}{ccc} y_{40} & y_{31} & y_{22} \end{array}
 \\

 \begin{array}{ccc} y_{11}\ | & y_{21} & y_{12}| \end{array}
\ \begin{array}{ccc} y_{31} & y_{22} & y_{13} \end{array}
 \\

 \begin{array}{ccc}y_{02}\ | & y_{12} & y_{03}| \end{array}
\ \begin{array}{ccc} y_{22} & y_{13} & y_{04} \end{array}

 \end{array}
\right]
\end{equation}
\textbf{Localizing matrix: }Given a polynomial $\mathcal{P} \in \mathbb R [x]$ with coefficient vector $ \mathbf p = \{ p_{\gamma }\}_{\gamma\in\mathbb{N}^n}$ and degree $\delta$, localizing matrix $M_r(\mathbf{py})$ with respect to $\mathbf y $ and $\mathbf p$ is as follows \cite{ref_Int_opt1, ref_Int_opt2}:
\begin{equation}\label{localization matrix def}
M_r(\mathcal P(x) \mathbf y )(i,j) =\sum_{\gamma \in \mathbb N^n} p_{\gamma} y_{\gamma +\alpha^{(i)}+\alpha^{(j)}}, \ \ 1 \leq i,j \leq  S_{n,r}.
\end{equation}
For example, given $\mathbf{y}=\{y_\alpha\}_{\alpha\in\mathbb{N}^2}$ and polynomial $\cP$,
\begin{equation}
  \mathcal{P}(x)=ax_1-bx^2_2,
\end{equation}
the localizing matrix for $r=1$ is formed as follows
\begin{equation}
  M_1( \mathcal{P}(x)\mathbf{y})= \begin{small}
    \left[ \begin{array}{ccc}
ay_{10}-by_{02} & ay_{20}-by_{12} & ay_{11}-by_{03} \\
ay_{20}-by_{12} & ay_{30}-by_{22} & ay_{21}-by_{13} \\
ay_{11}-by_{03} & ay_{21}-by_{13} & ay_{12}-by_{04} \end{array}
\right]
  \end{small}
\end{equation}

\subsection{ Preliminary Results}

In this section, we state some standard results found in the literature that will be referred to later. The following results give necessary, sufficient conditions for $\mathbf y$ to have a representing measure $\mu$ -- for details see \cite{ref_PreR_2,ref_Int_opt1, ref_Int_opt2}.

Consider the semialgebraic set $\cK$ defined as
\begin{equation}\label{preliminary result_semi algebraic set}
\cK = \{ x\in \mathbb{R}^n: g_j(x)\geq0,\ j=1,2,\dots ,\ell\ \}.
\end{equation}
for some polynomials $\mathcal{P}_j\in \mathbb R [x]$, and assume that $\cK$ satisfies Putinar's property.

(i) If $f \in \mathbb R [x]$ is strictly positive on $\cK$, then:
\begin{equation} 
f = \sigma_0 + \sum_{j=1}^{l} \sigma_jg_j
\end{equation} 
for some $\sigma_j \in \sum^2[x]$.\\

(ii) The sequence $\mathbf y = \{y_{\alpha}\}_{\alpha \in \mathbb{N}^n}$ has a representing finite Borel measure $\mu $ on $\cK$ if and only if:

\begin{equation} 
\label{PR_moment}
M_r(\mathbf y)\succcurlyeq 0, M_r(g_j \mathbf y )\succcurlyeq 0,\ \ j=1,\dots ,m 
\end{equation} 
for every $ r\in \mathbb N^n$. \\


\section{ Convex Formulation} \label{Sec_Convex}

The approach presented in this paper relies on finding polynomial approximations of the indicator function of the support set of the measure of interest. More precisely, 
let $\cK$ represent the support set of a given measure $\mu$. 
The results in this paper aim at finding polynomial approximations of
\begin{center}
$\mathbb{I}_{\cK}(x) \doteq \left\{ \begin{array}{cc}
1 & \mbox{if $x \in \cK$}\\
0 & \mbox{otherwise}.\end{array} \right.  $
\end{center}
and use the level sets of these polynomials to approximate~$\cK$. 
In order to approximate the indicator function above consider the following optimization problem. \\

\begin{Problem} \label{Prob2}
Let $d$ be a given integer. Moreover, let $\mathcal{B}$ be a known (simple) set containing the support set $\cK$ and  $\mu_{\mathcal{B}}$ be the Lebesgue measure supported on the set $\mathcal{B}$. Solve
\begin{align}
\label{Convex1}
\mathbf{P_2^*}:=&\ \min_{\mathcal{P}_d(x) \in \mathbb R_{\rm d}[x] } \int {\mathcal{P}_d(x)} d\mu_{\mathcal{B}} \\
\hbox{s.t.}\quad & \mathcal{P}_d(x)\geq 0, x \in \cB \subeqn\\
& \mathcal{P}_d(x)\geq 1, x \in \cK \subeqn.
\end{align}
\end{Problem}

For every $d$, the problem above provides a polynomial $\mathcal{P}_d^*$ with the smallest $\ell_1$-norm on $\cB$ that is i) positive in the (simple) set $\cB$ and ii) larger than one  inside the support set $\cK$. For this (infinite dimensional) optimization problem we have the following result.

\vskip .1in
\begin{Theo}
 For a given integer $d$, let 
\[ \mathcal{K}_d \doteq \{ x \in \mathbb{R}^n: \mathcal{P}_d^*(x) \geq 1 \}\] 
 be the semialgebraic set constructed using the solution $ \mathcal{P}_d^*$ of the problem (\ref{Convex1}). Then 
\[ \lim_{\rm d \rightarrow \infty } \mu_{\mathcal{B}}(\mathcal{K}_d - \mathcal{K}) = 0. \]
\end{Theo}

\vskip .1in
\begin{proof} %
As in \cite{ref_Indicator} one can show that $\cP^*_d$ converges almost uniformly (with respect to measure $\mu_{\mathcal{B}}$) to the indicator function $\mathbb{I}_{\cK}$. Moreover, one has $\cK \subseteq \cK_d$ for all $d$.
These two facts imply that 
\[ \lim_{\rm d \rightarrow \infty } \mu_{\mathcal{B}}(\mathcal{K}_d - \mathcal{K}) = 0 \]
which completes the proof.
\end{proof}

In the optimization problem above, one approximates the indicator function of the set $\cK$ by using the knowledge that  this set is contained in a known set $\cB$. This set is usually chosen in such a way that one can compute all the moments of the measure $\mu_{\mathcal{B}}$ in a closed form.

However, the problem above obviously requires the knowledge of the measure $\mu$ whose support $\cK$ we are trying to determine. To be able to solve this problem by using only knowledge of moments consider a bounding  set $\mathcal{B}$ defined by a set of polynomial inequalities; i.e.,
\[\mathcal{B} = \left\lbrace   x \in \mathbb{R}^n: g_j(x) \geq 0, j=1,...,l \right\rbrace  \]
where $g_j$, $j=1,2,\ldots,l$ are given polynomials. As before, let $ \mu_{\mathcal{B}}$ be the Lebesgue measure supported in $\cB$ with \mbox{$\alpha$-th} moment $y_{\mathcal{B}_{\alpha}}$.  Moreover, let the (infinite) vector $\textbf{y}$ be the vector containing all the moments of the measure $\mu$. Then, define the following optimization problem (which has an infinite number of constraints).

\begin{Problem}  \label{Prob3}
\begin{align}
\label{SDP1}
\mathbf{P_3^*}:=&\ \min_{p_{\alpha}, \sigma_j} \sum_{\alpha=0}^d p_{\alpha} y_{\mathcal{B}_{\alpha}} \\
\hbox{s.t.}\quad & \cP_d(x) = \sum_{\|\alpha\|_1 \leq d} p_\alpha x^\alpha \subeqn\\
& \mathcal{P}_d(x) = \sigma_0(x) + \sum_{j=1}^{l} \sigma_j(x)g_j(x)   \subeqn \label{subeq:p3_sos}\\
& \sigma_j \in \Sigma^2[x]; j=0,1,\ldots,l \subeqn\\
& deg(\sigma_0) \leq d; \subeqn\\
& deg(\sigma_j g_j) \leq d;  j=1,2,\ldots,l \subeqn\\
& M_{\infty}((\mathcal{P}_d(x)-1)\textbf{y})\succcurlyeq 0  \label{subeq:p3_localization}\subeqn
\end{align}
\end{Problem}

The problem above is a first step towards an implementable version of Problem \ref{Prob2}. The objective function is the same in both, just represented as a function of the moments of $\mu_\cB$ in Problem \ref{Prob3}. Constraint~\eqref{subeq:p3_sos} enforces $\cP_d$ to be positive on the set $\cB$. Finally, given the definition of localization matrix, constraint~\eqref{subeq:p3_localization} ensures that $\cP_d$ is larger than one in the support set of $\mu$.

Since one cannot solve the problem above, in this paper we propose the following relaxation.

\begin{Problem} \label{Prob4}
\begin{align}
\label{SDP2}
\mathbf{P_4^*}:=&\ \min_{p_{\alpha}, \sigma_j} \sum_{\alpha=0}^d p_{\alpha} y_{\mathcal{B}_{\alpha}} \\
\hbox{s.t.}\quad & \cP_d(x) = \sum_{\|\alpha\|_1 \leq d} p_\alpha x^\alpha \subeqn\\
& \mathcal{P}_d(x) = \sigma_0(x) + \sum_{j=1}^{l} \sigma_j(x)g_j(x) \subeqn \label{subeq:p4sos}\\
& \sigma_j \in \Sigma^2[x]; j=0,1,\ldots,l \subeqn\\
& deg(\sigma_0) \leq d; \subeqn\\
& deg(\sigma_j g_j) \leq d;  j=1,2,\ldots,l \subeqn\\
& M_{r}((\mathcal{P}_d(x)-1)\textbf{y})\succcurlyeq 0  \label{subeq:4_localization}\subeqn
\end{align}
\end{Problem}

\begin{figure}[!h]
 \centering
 \includegraphics[keepaspectratio=true,scale=0.3]{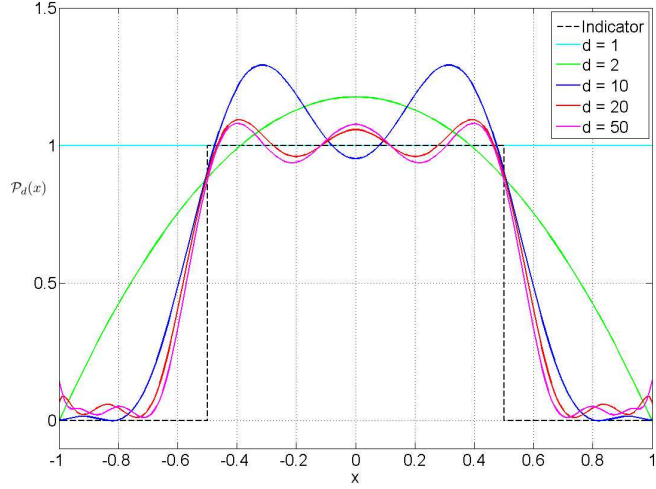}
 \caption{Result of SDP in (\ref{SDP2}) For Example \ref{Exa1}} \label{Fig_Exa1}
\end{figure}

where, $r \geq$1 is relaxation order. In other words, we truncate the infinite moment localization matrix. One should note that the problem above can be formulated as a standard SDP; i.e., minimization of a linear function subject to Linear Matrix Inequalities (LMIs).

The truncation of the moment localization matrix  provides an approximation of the constraint $ \mathcal{P}_d(x)\geq 1$ for all $x \in \cK $. Although, if $r$ is ``large'' one has acceptable estimates of the support set, for ``low'' values of $r$ this can lead to estimates of the support set that are less accurate than desirable.
\\

\begin{Exa} \label{Exa1}
Let, $\textbf{y}$ be a moment sequence of uniform probability measure $\mu$ supported on $[-0.5, 0.5]$. The $\alpha$-th moment of uniform distribution $U[a,b]$ is $y_{\alpha} =\frac{b^{\alpha+1} - a^{\alpha+1}}{(b-a)(\alpha+1)} $. For this example, we take $\mathcal{B}=[-1, 1]$, and use the moments up to order 2$d$. To solve the SDP (\ref{SDP2}), Yalmip is used which is a Matlab-based toolbox aimed at optimization \cite{ref_Yalmip}. The obtained results are depicted in Fig \ref{Fig_Exa1}. One can see as $d$, the order of polynomial, increases $\mathcal{P}_d(x)$ converges to indicator function of support of uniform measure. Hence, the semialgebraic set $\mathcal{K}_d = \{ x \in \mathbb{R}: \mathcal{P}_d(x) \geq 1 \}$ provides better approximations of the support as one increases $d$. However, as one can see in Fig \ref{Fig_Exa1}, $\cP_d$ can be below one in a significant subset of the support of $\mu$.
\end{Exa}

\subsection{An Heuristic for Improved Performance}
To minimize the measure of the subset of the support of the measure $\mu$ where $\cP_d$ is below one, we propose to maximize the values of $\mathcal{P}_d(x)$ inside the support of the measure while still trying to bring its values as low as possible everywhere else in $\cB$. This results in following modified SDP.\\

\begin{Problem} \label{Prob5}
\begin{align}
\label{SDP3}
\mathbf{P_5^*}:=&\ \min_{p_{\alpha}, \sigma_j} \sum_{\alpha=0}^d p_{\alpha} y_{\mathcal{B}_{\alpha}} - \omega_h h \\
\hbox{s.t.}\quad & \cP_d(x) = \sum_{\|\alpha\|_1 \leq d} p_\alpha x^\alpha \subeqn\\
& \mathcal{P}_d(x) = \sigma_0(x) + \sum_{j=1}^{l} \sigma_j(x)g_j(x)  \subeqn \label{subeq:p4sos}\\
& \sigma_j \in \Sigma^2[x]; j=0,1,\ldots,l \subeqn\\
& deg(\sigma_0) \leq d; \subeqn\\
& deg(\sigma_j g_j) \leq d;  j=1,2,\ldots,l \subeqn\\
& M_{r}((\mathcal{P}_d(x)-h)\textbf{y})\succcurlyeq 0 \label{SDP1_3} \subeqn\\
& 1 \leq h \leq 1+ \Delta h \subeqn
\end{align}
\end{Problem}
where, $\omega_h$ and $\Delta h$ are positive design parameters.

\begin{figure}
 \centering
 \includegraphics[keepaspectratio=true,scale=0.3]{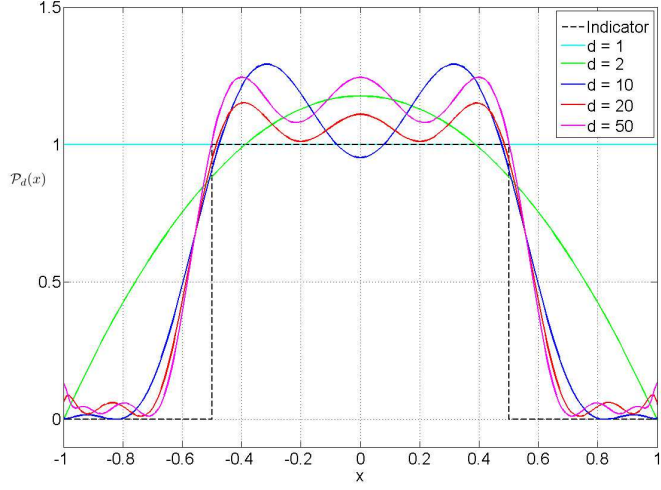}
 \caption{Result of SDP in (\ref{SDP3}) For Example \ref{Exa1}} \label{Fig_Exa2}
\end{figure}

 To show the effectiveness of the modified SDP, we again consider the uniform measure in Example \ref{Exa1}. Fig \ref{Fig_Exa2} shows the results obtained by solving the modified SDP with parameters $\omega_h = 1.2$ and $\Delta h = 0.2 $. As it is seen, on obtains a substantial improvement in the estimate of the support set. \\

\begin{figure}[!h]
 \centering
 \includegraphics[keepaspectratio=true,scale=0.3]{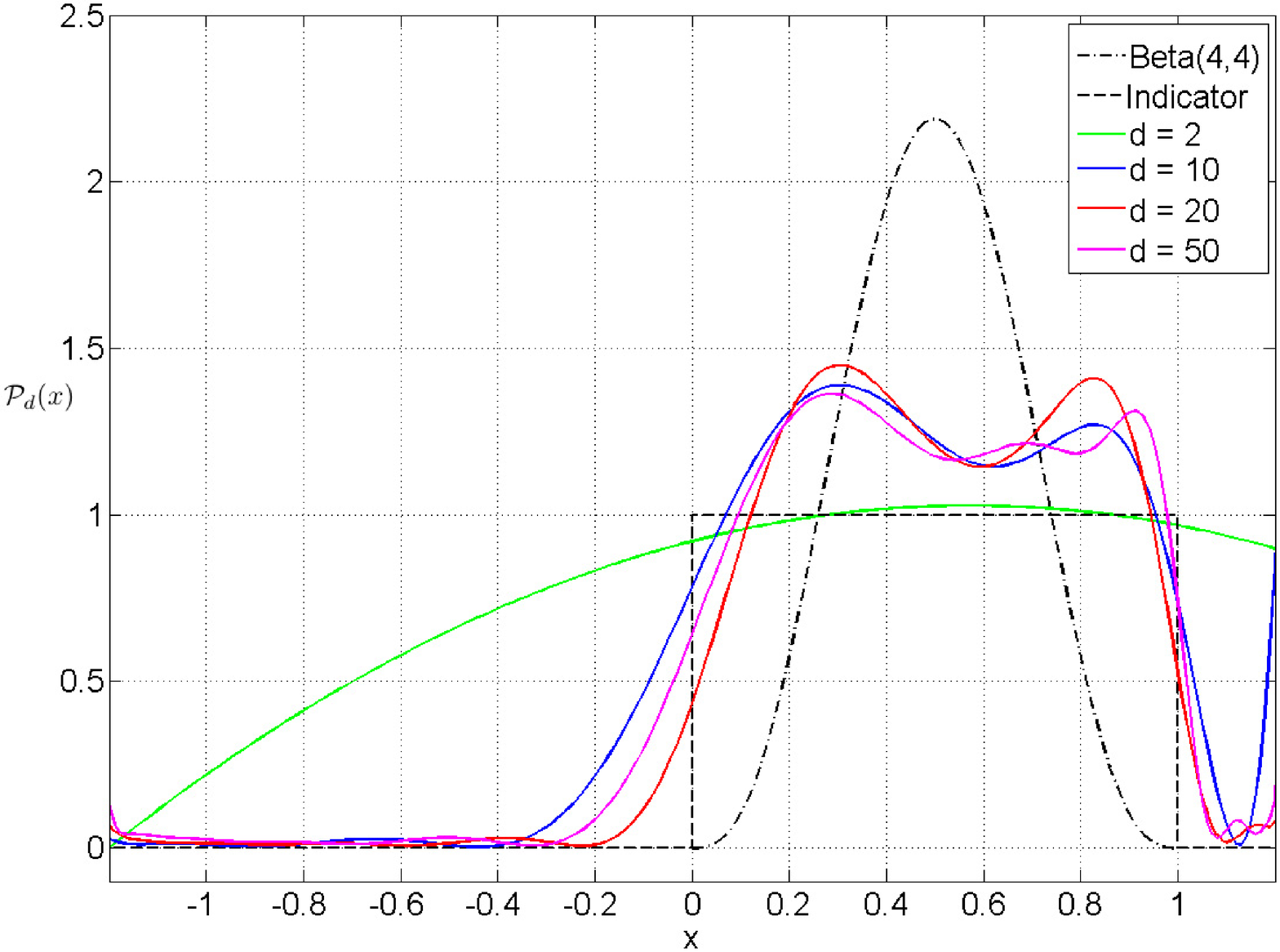}
 \caption{Result of SDP in (\ref{SDP3}) For Example \ref{Exa2}} \label{Fig_Exa3}
\end{figure}

\begin{Exa}\label{Exa2}
In this example, we consider a  $\mathrm{Beta}(4,4)$  probability measure on $[0, 1]$. The $\alpha$-th moment of Beta distribution $\mathrm{Beta}(a,b)$ over [0,1] is $y_{\alpha} = \frac{a+k-1}{(a+b+\alpha-1)} y_{\alpha-1}$ and $y_0=1$. We assume that set $\mathcal{B}=[-1.2, 1.2]$, and use the moments up to order 2$d$. The obtained results by solving SDP  (\ref{SDP3}) with parameters $\omega_h = 1.2$ and $\Delta h = 0.2 $ are depicted in Fig \ref{Fig_Exa3}. 
\end{Exa}

This is a more difficult problem than previous ones since, in terms of probability, there is a ``smooth transition''  from the interior to the exterior of the support set. Nevertheless, if one uses enough moments, one can get a very good approximation of the support.\\

\begin{figure}[!h] 
 \centering
 \includegraphics[keepaspectratio=true,scale=0.29]{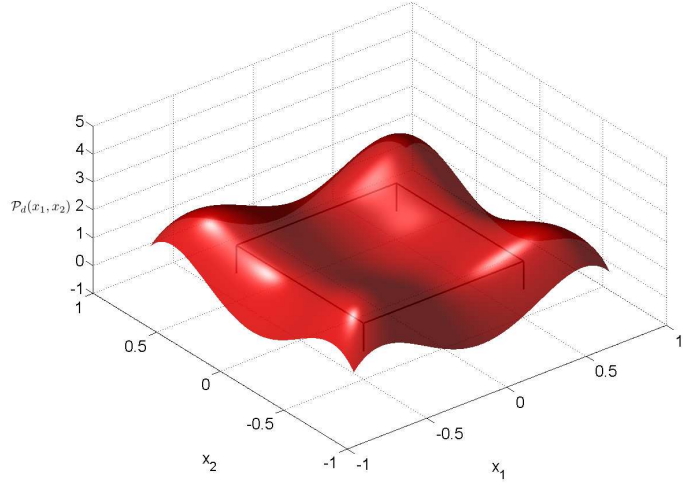}
 \caption{Result of SDP in (\ref{SDP3}) For Example \ref{Exa3} } \label{Fig_Exa4}
\end{figure}

\begin{Exa}\label{Exa3}
In here, we consider a 2-dimensional example where one wants to approximate the support of a  uniform probability measure on $[-0.5, 0.5]^2$. The results obtained by solving SDP (\ref{SDP3}) with parameters $d = 14 $, $\omega_h = 1.2$, and $\Delta h = 0.2 $ are depicted in Fig \ref{Fig_Exa4}.
\end{Exa}

\section{Support Reconstruction for Uniform Measures} \label{Sec_Uniform}

In this section, we present a modification of our approach aimed specifically at uniform distributions. In the development to follow, we rely on a result in \cite{ref_Pre_3} which provides criteria under which  polynomials vanish on the boundary of support of the uniform measure of interest. We now elaborate on this.

Define  
\begin{equation}\label{momnt matirx def}
\bar{M}_r ( \mathbf y )(i,j)=\dfrac{n+|i|+|j|}{n+|i|} y_{\alpha^{(i)}+\alpha^{(j)}},\ \ \ 1 \leq i,j \leq S_{n,r},
\end{equation}
where $ \mathbf y =\{ y_\alpha\}_{\alpha\in\mathbb{N}^n}$ are the moments  of the uniform distribution of interest. The results in \cite{ref_Pre_3} show that a  polynomial $\mathcal{P}(x)$ whose vector of coefficients $\textbf{p}$ is the eigenvector associated with zero eigenvalue of the matrix $\bar{M}_r$, vanishes on the boundary of support of measure. More precisely, under some technical conditions, 
\begin{equation}
\bar{M}_r(\textbf{y}) \textbf{p} = 0 \Rightarrow \cP(x) = 0 \text{ for all } x \in \partial \mathcal{K}
\end{equation}
where, $\partial \mathcal{K}$ denotes the boundary of support set $\mathcal{K}$. However, without any additional constraints, this polynomial can also be zero in the interior of $\cK$ and, hence, it might not provide a good estimate of the support.

Nevertheless, one can  take advantage of this property and modify our approach as follows. \\

\begin{Problem} \label{Prob6}
\begin{align}
\label{SDP4}
\mathbf{P_6^*}:=&\ \min_{p_{\alpha}, \sigma_j} \sum_{\alpha=0}^d p_{\alpha} y_{\mathcal{B}_{\alpha}} - \omega_h h + \omega_M \Vert \bar{M}_d(\textbf{y}) (\textbf{p}-1) \Vert_2 \\
\hbox{s.t.}\quad & \cP_d(x) = \sum_{\|\alpha\|_1 \leq d} p_\alpha x^\alpha \subeqn\\
& \mathcal{P}_d(x) =\sigma_0(x) + \sum_{j=1}^{l} \sigma_j(x)g_j(x) \subeqn \label{subeq:p5sos}\\
& \sigma_j \in \Sigma^2[x]; j=0,1,\ldots,l \subeqn\\
& deg(\sigma_0) \leq d; \subeqn\\
& deg(\sigma_j g_j) \leq d;  j=1,2,\ldots,l \subeqn\\
& M_{r}((\mathcal{P}_d(x)-h)\textbf{y})\succcurlyeq 0 \label{SDP1_3} \subeqn\\
& 1 \leq h \leq 1+ \Delta h \subeqn
\end{align}
\end{Problem}
where, $\omega_M$, $\omega_h$ and $\Delta h$ are positive design parameters, \mbox{$\mathbf{p}=\{p_\alpha\}_{\alpha\in\mathbb{N}^n}$} denotes the vector of polynomial coefficients and $ \Vert. \Vert_2$ denotes the $l_2$ norm. 

In fact in (\ref{SDP4}), we  aim at ``pushing'' the coefficients of the  polynomial $(\mathcal{P}(x)-1)$ as close as possible to the null space of  $\bar{M}_d$ by minimizing the term $\Vert \bar{M}_d (\textbf{p}-1) \Vert_2$. In this case obtained polynomial $\mathcal{P}_d(x)$ becomes close to one at the boundary of support while we still aim at having $\cP_d
$ larger than one inside the support. 

\begin{figure}
 \centering
 \includegraphics[keepaspectratio=true,scale=0.3]{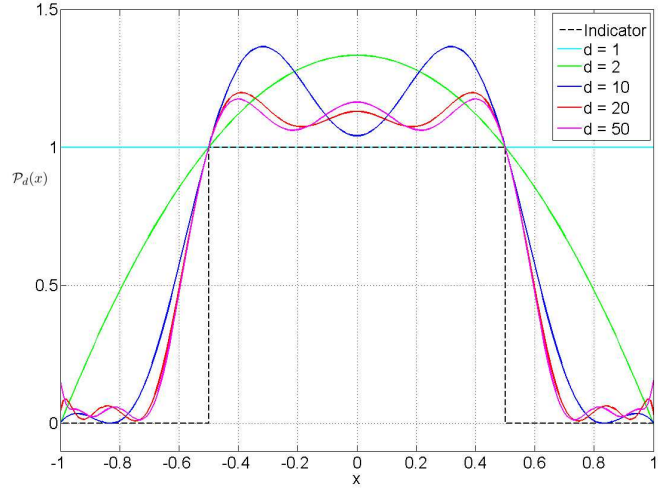}
 \caption{Result of SDP in (\ref{SDP4}) For Example \ref{Exa1}} \label{Fig_Exa5}
\end{figure}

To show the effectiveness of proposed method, we reconstruct the support for the measure of Example \ref{Exa1} by solving the SDP  (\ref{SDP4}) with parameter $\omega_M = 10$. The obtained result are depicted in Fig \ref{Fig_Exa5}, where semialgebraic set $\mathcal{K}_d = \{ x \in \mathbb{R}: \mathcal{P}(x)_d \geq 1\}$ for any polynomial order $d \geq 2$ exactly reconstructs the support of measure.\\

\begin{figure}[!h]
 \centering
 \includegraphics[keepaspectratio=true,scale=0.3]{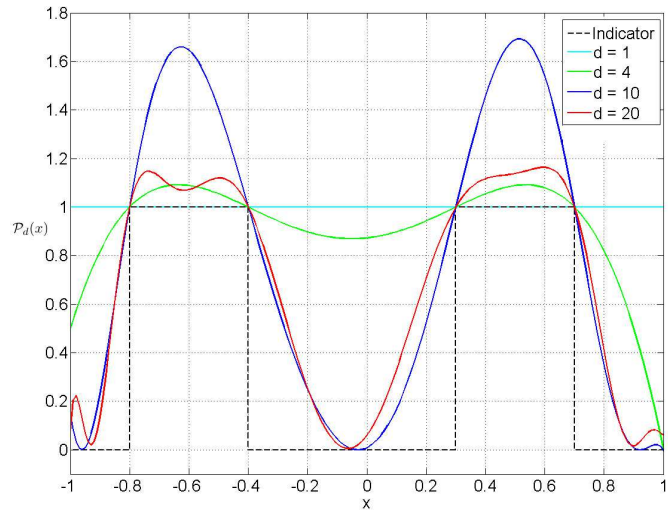}
 \caption{Result of SDP in (\ref{SDP4}) For Example \ref{Exa4} } \label{Fig_Exa6}
\end{figure}

\begin{Exa}\label{Exa4}
To further show the effectiveness of our approach, we now consider a uniform distribution with disconnected support. More precisely, we aim at estimating the support of a uniform probability measure over the union of the sets $[-0.8, -0.4]$ and $[0.3, 0.7]$. We assume that $\mathcal{B}=[-1, 1]$ and use  moments up to order $2d$. The results obtained by solving SDP  (\ref{SDP4}) with parameters $\omega_h = 1.2$, $\omega_M = 10$ and $\Delta h = 0.2 $ are depicted in Fig \ref{Fig_Exa6}, where one can see that the semialgebraic set $\mathcal{K}_d = \{ x \in \mathbb{R}: \mathcal{P}_d(x) \geq 1\}$ for   $d \geq 4$ exactly reconstructs the support of measure.
\end{Exa}

\section{ Conclusion} \label{Sec_Con}

In this paper, we present a novel approach to the problem of reconstruction of support of measures from their moments. A sequence of semidefinite relaxations is provided whose solution converge to the support of the measure of interest. Examples are provided that show that one does obtain a good approximation of support using only a finite number of moments. Further research effort is now being put  on developing methods for support reconstruction for specific classes of measures which have provable performance.


\bibliography{MyRef}
\bibliographystyle{abbrv}

\end{document}

%% file: defs.tex
\def\cB{\mathcal{B}}

\def\cK{\mathcal{K}}

\def\cM{\mathcal{M}}

\def\cP{\mathcal{P}}

\def\smskip{\smallskip}

\def\texitem#1{\par\smskip\noindent\hangindent 25pt
               \hbox to 25pt {\hss #1 ~}\ignorespaces}

\def\norm#1{\|#1\|}

\newcommand{\BEAS}{\begin{eqnarray*}}
\newcommand{\EEAS}{\end{eqnarray*}}
\newcommand{\BEA}{\begin{eqnarray}}
\newcommand{\EEA}{\end{eqnarray}}
\newcommand{\BEQ}{\begin{eqnarray}}
\newcommand{\EEQ}{\end{eqnarray}}
\newcommand{\BIT}{\begin{itemize}}
\newcommand{\EIT}{\end{itemize}}
\newcommand{\BNUM}{\begin{enumerate}}
\newcommand{\ENUM}{\end{enumerate}}

\newcommand{\BA}{\begin{array}}
\newcommand{\EA}{\end{array}}

\newcommand{\reals}{\mathbb{R}}

\newif\ifpagenumbering
\pagenumberingtrue

\pagenumberingfalse

\newsavebox{\theorembox}
\newsavebox{\lemmabox}
\newsavebox{\remarkbox}
\newsavebox{\assbox}
\savebox{\theorembox}{\noindent\bf Theorem}
\savebox{\lemmabox}{\noindent\bf Lemma}
\savebox{\remarkbox}{\noindent\bf Remark}
\savebox{\assbox}{\noindent\bf Assumption}

